
%
%
%
%
%
%
%


\input amssym.def
\input amssym

\magnification=1100
\overfullrule=0pt

\font\eightrm=cmr8

\textfont0=\tenrm
\scriptfont0=\sevenrm
\scriptscriptfont0=\fiverm

%
%
\font\teni=cmmi10
\font\seveni=cmmi7
\font\fivei=cmmi5
\textfont1=\teni
\scriptfont1=\seveni
\scriptscriptfont1=\fivei

%

%

\font\eightit=cmti8

%
\font\eightbfti=cmbxti10 at 8pt
\font\bfti=cmbxti10


\font\sc=cmcsc10


%
\font\tenbb=msbm10
\font\sevenbb=msbm7
\font\fivebb=msbm5
\newfam\bbfam
\textfont\bbfam=\tenbb
\scriptfont\bbfam=\sevenbb
\scriptscriptfont\bbfam=\fivebb
\def\bb{\fam\bbfam\tenbb}
\let\oldbb=\bb
\def\bb #1{{\oldbb #1}}

%

%
\font\tengoth=eufm10
\font\sevengoth=eufm7
\font\fivegoth=eufm5
\newfam\gothfam
\textfont\gothfam=\tengoth
\scriptfont\gothfam=\sevengoth
\scriptscriptfont\gothfam=\fivegoth

%

%

\font\fourteenbf=cmbx10 at 14pt

\font\eightbf=cmbx8
\font\sixbf=cmbx5 at 6pt 

\font\tenbf=cmbx10
\font\sevenbf=cmbx7
\font\fivebf=cmbx5
\newfam\bffam
\textfont\bffam=\tenbf
\scriptfont\bffam=\sevenbf
\scriptscriptfont\bffam=\fivebf
\def\bf{\fam\bffam\tenbf}
%
%



\newif\ifpagetitre           \pagetitretrue
\newtoks\hautpagetitre       \hautpagetitre={\hfil}
\newtoks\baspagetitre        \baspagetitre={\hfil}

\newtoks\auteurcourant       \auteurcourant={\hfil}
\newtoks\titrecourant        \titrecourant={\hfil}

\newtoks\hautpagegauche
         \hautpagegauche={\hfil\the\auteurcourant\hfil}
\newtoks\hautpagedroite
         \hautpagedroite={\hfil\the\titrecourant\hfil}

\newtoks\baspagegauche
         \baspagegauche={\hfil\tenrm\folio\hfil}
\newtoks\baspagedroite
         \baspagedroite={\hfil\tenrm\folio\hfil}

\headline={\ifpagetitre\the\hautpagetitre
\else\ifodd\pageno\the\hautpagedroite
\else\the\hautpagegauche\fi\fi }

\footline={\ifpagetitre\the\baspagetitre
\global\pagetitrefalse
\else\ifodd\pageno\the\baspagedroite
\else\the\baspagegauche\fi\fi }



\def\og{\leavevmode\raise.3ex
     \hbox{$\scriptscriptstyle\langle\!\langle$~}}
\def\fg{\leavevmode\raise.3ex
     \hbox{~$\!\scriptscriptstyle\,\rangle\!\rangle$}}
     


\def\build #1_#2^#3{\mathrel{\mathop{\kern 0pt #1}
     \limits_{#2}^{#3}}}


\def\scs{\scriptstyle}



\def\hfl[#1][#2][#3]#4#5{\kern -#1
 \sumdim=#2 \advance\sumdim by #1 \advance\sumdim by #3
  \smash{\mathop{\hbox to \sumdim {\rightarrowfill}}
   \limits^{\scriptstyle#4}_{\scriptstyle#5}}
    \kern-#3}

\def\vfl[#1][#2][#3]#4#5%
 {\sumdim=#1 \advance\sumdim by #2 \advance\sumdim by #3
  \setbox1=\hbox{$\left\downarrow\vbox to .5\sumdim {}\right.$}
   \setbox1=\hbox{\llap{$\scriptstyle #4$}\box1
    \rlap{$\scriptstyle #5$}}
     \vcenter{\kern -#1\box1\kern -#3}}


\newdimen\sumdim
\def\diagram#1\enddiagram
    {\vcenter{\offinterlineskip
      \def\tvi{\vrule height 10pt depth 10pt width 0pt}
       \halign{&\tvi\kern 5pt\hfil$\displaystyle##$\hfil\kern 5pt
        \crcr #1\crcr}}}







\def\vspace[#1]{\noalign{\vskip #1}}



\def\boxit [#1]#2{\vbox{\hrule\hbox{\vrule
     \vbox spread #1{\vss\hbox spread #1{\hss #2\hss}\vss}%
      \vrule}\hrule}}



\catcode`\@=11
\def\Eqalign #1{\null\,\vcenter{\openup\jot\m@th\ialign{
\strut\hfil$\displaystyle{##}$&$\displaystyle{{}##}$\hfil
&&\quad\strut\hfil$\displaystyle{##}$&$\displaystyle{{}##}$
\hfil\crcr #1\crcr}}\,}
\catcode`\@=12






\vsize= 20.5cm
\hsize=15cm
\hoffset=3mm
\voffset=3mm


\def\tend_#1^#2{\mathrel{
	\mathop{\kern 0pt\hbox to 1cm{\rightarrowfill}}
	\limits_{#1}^{#2}}}

\def\lfq{\leavevmode\raise.3ex\hbox{$\scriptscriptstyle
\langle\!\langle$}\thinspace}
\def\rfq{\leavevmode\thinspace\raise.3ex\hbox{$\scriptscriptstyle
\rangle\!\rangle$}}

\font\msamten=msam10
\font\msamseven=msam7
\newfam\msamfam

\textfont\msamfam=\msamten
\scriptfont\msamfam=\msamseven

\def\hexnumber #1%
{\ifcase #1 0\or 1\or 2\or 3\or 4\or 5
   \or 6\or 7\or 8\or 9\or A\or B\or C\or D\or E\or F
    \fi}

\mathchardef\leadsto="3\hexnumber\msamfam 20

\def\frac #1#2{{#1\over #2}}

\def\qed{\quad\raise -2pt\hbox{\vrule\vbox to 10pt{\hrule width 4pt 
\vfill\hrule}\vrule}} 

\def\qeda{\quad\raise -2pt\hbox{\vrule\vbox to 10pt{\hrule 
width 4pt height9pt 
\vfill\hrule}\vrule}}




\auteurcourant{\eightbf B.~Lass \hfill}


\titrecourant{\hfill\eightbf D\'emonstration de la conjecture de Dumont}


{\parindent=0pt\eightbf
C. R. Acad. Sci. Paris, t.~xxx, S\'erie~I,
p.~xxx--xxx, 2005


Th\'eorie des nombres/\hskip -.5mm{\eightbfti Number Theory}

(Combinatoire/\hskip -.5mm{\eightbfti Combinatorics})

\vskip 1.5cm



{\fourteenbf
D\'emonstration de la conjecture de Dumont
}
\parindent=-1mm\footnote{}{
Note pr\'esent\'ee par \'Etienne G{\sixbf HYS.}
}

\vskip 15pt


\noindent{\bf Bodo LASS}

\vskip 6pt


{\parindent=0mm
Institut Camille Jordan, UMR 5208 du CNRS, Universit\'e Claude Bernard Lyon 1, 

43, Boulevard du 11 Novembre 1918, 69622 Villeurbanne cedex, France

Courriel~: lass@igd.univ-lyon1.fr
\par}

\medskip


(Re\c{c}u le 8 juillet 2005, accept\'e le xx xxxx 2005)
}

\medskip

\centerline{\eightbfti \`A Dominique Foata, pour son 70-i\`eme anniversaire}

\vskip3pt

\line{\hbox to 2cm{}\hrulefill\hbox to 1.5cm{}}

\vskip3pt


{\parindent=2cm\rightskip 1.5cm\baselineskip=9pt
\item{\bf R\'esum\'e.~}{\eightrm Soit 
\vskip-10pt
$$\openup-5pt\eqalign{
\scs r_k^{1(2)}(n) \; := \; &\scs|\{(x_1,x_2,\dots,x_k) \in {\bb N}^k \; | \;
     n \; = \; x_1^2+x_2^2+\cdots+x_k^2, \; {\scs x_i \equiv 1 \!\!\! \pmod 2}, \; 1 \le i \le k\}|, \cr 
\scs c_k^{1(4)}(n) \; := \; &\scs|\{(x_1,x_2,\dots,x_k) \in {\bb N}^k \; | \;
     n \; = \; x_1x_2+x_2x_3+\cdots+x_{k-1}x_{k}+x_{k}x_1, \; {\scs x_i \equiv 1 \;(4)}\}|, \cr 
\scs c_k^{3(4)}(n) \; := \; &\scs|\{(x_1,x_2,\dots,x_k) \in {\bb N}^k \; | \;
     n \; = \; x_1x_2+x_2x_3+\cdots+x_{k-1}x_{k}+x_{k}x_1, \; {\scs x_i \equiv 3 \;(4)}\}|. \cr 
}$$
\vskip-5pt
Dumont~[2] a conjectur\'e l'identit\'e $\scs r_k^{1(2)}(n) = c_k^{1(4)}(n) - (-1)^k c_k^{3(4)}(n)$
qui g\'en\'eralise, notamment, les r\'esultats classiques de Lagrange, Gau\ss, Jacobi et Kronecker 
sur les d\'ecomposi\-tions de tout entier en deux, trois et quatre carr\'es. Nous donnons une preuve 
combinatoire de la conjecture de Dumont.~\copyright~Acad\'emie des Sciences, Paris}
\par}

\vskip 10pt


{\parindent=2cm\rightskip 1.5cm
{\bfti A proof of Dumont's conjecture}
\par}

\vskip 10pt


{\parindent=2cm\rightskip 1.5cm\baselineskip=9pt
\item{\bf Abstract.~}{\eightit Let 
\vskip-10pt
$$\openup-5pt\eqalign{
\scs r_k^{1(2)}(n) \; := \; &\scs|\{(x_1,x_2,\dots,x_k) \in {\bb N}^k \; | \;
     n \; = \; x_1^2+x_2^2+\cdots+x_k^2, \; {\scs x_i \equiv 1 \!\!\! \pmod 2}, \; 1 \le i \le k\}|, \cr 
\scs c_k^{1(4)}(n) \; := \; &\scs|\{(x_1,x_2,\dots,x_k) \in {\bb N}^k \; | \;
     n \; = \; x_1x_2+x_2x_3+\cdots+x_{k-1}x_{k}+x_{k}x_1, \; {\scs x_i \equiv 1 \;(4)}\}|, \cr 
\scs c_k^{3(4)}(n) \; := \; &\scs|\{(x_1,x_2,\dots,x_k) \in {\bb N}^k \; | \;
     n \; = \; x_1x_2+x_2x_3+\cdots+x_{k-1}x_{k}+x_{k}x_1, \; {\scs x_i \equiv 3 \;(4)}\}|. \cr 
}$$
\vskip-5pt
Dumont~[2] has conjectured the identity $\scs r_k^{1(2)}(n) = c_k^{1(4)}(n) - (-1)^k c_k^{3(4)}(n)$,
which generalizes, in particular, the classical results of Lagrange, Gau\ss, Jacobi and Kronecker 
on the sums of two, three and four squares. We give a combinatorial proof of Dumont's 
conjecture.~}{\eightrm\copyright~Acad\'emie des
Sciences, Paris}
\par}

\vskip1.5pt
\line{\hbox to 2cm{}\hrulefill\hbox to 1.5cm{}}
\medskip

\parindent=3mm
\vskip 5mm


\noindent
{\bf Abridged English Version}

\bigskip

\noindent
Let ${\bb N}^{a(b)}$ denote the set of strictly positive integers congruent to $a$ modulo $b$,
and let $k$ and $n$ be integers such that $1 \le k \le n$ and $n \equiv k \!\!\! \pmod 8$.
We want to count the number of solutions of the following Diophantine equations:
$$\openup2pt\eqalign{
r_k^{1(2)}(n) \; &:= \; |\{(x_1,x_2,\dots,x_k) \in ({\bb N}^{1(2)})^k \; | \;
     n \; = \; x_1^2+x_2^2+\cdots+x_k^2 \}|, \cr 
c_k^{1(4)}(n) \; &:= \; |\{(x_1,x_2,\dots,x_k) \in ({\bb N}^{1(4)})^k \; | \;
     n \; = \; x_1x_2+x_2x_3+\cdots+x_{k-1}x_{k}+x_{k}x_1 \}|, \cr 
c_k^{3(4)}(n) \; &:= \; |\{(x_1,x_2,\dots,x_k) \in ({\bb N}^{3(4)})^k \; | \;
     n \; = \; x_1x_2+x_2x_3+\cdots+x_{k-1}x_{k}+x_{k}x_1 \}|. \cr 
}$$
Let $y_1,y_2,\dots,y_k$ be arbitrary strictly positive integers. The equivalence
$$
m \; = \; {y_1 \choose 2} + {y_2 \choose 2} + \cdots + {y_k \choose 2} 
\quad \Leftrightarrow \quad 
8m+k \; = \; (2y_1-1)^2 + (2y_2-1)^2 + \cdots + (2y_k-1)^2 
$$
shows that the congruence $n \equiv k \!\!\! \pmod 8$ is satisfied automatically in the first 
Diophantine equation (this is true for the two other equations too) and that $r_k^{1(2)}(n)$ 
also counts the number of decompositions of $(n-k)/8$ into $k$ triangular numbers. 
Moreover, Jacobi's two and four odd squares theorems tell us that for $n \equiv 2 \!\!\! \pmod 8$ 
and $n \equiv 4 \!\!\! \pmod 8$ we have 
$$
r_2^{1(2)}(n) = \sum_{d|(n/2)} (-1)^{(d-1)/2}, \qquad 
r_4^{1(2)}(n) = \sum_{d|(n/4)} d,
$$
respectively (the sums go over all positive divisors), whereas Kronecker's three odd squares 
theorem gives
$$
r_3^{1(2)}(n) = |\{(a,b,c) \in {\bb N}^3 \; | \; n = 4ac-b^2, b > 0, b < 2a, b < 2c\}|
$$
for $n \equiv 3 \!\!\! \pmod 8$. 
The aim of Andr\'e Weil's article~[3] is to prove exactly those three theorems.
Up to now, we could not see in such results anything else than special cases
of the general conjecture that number theory is less beautiful than
combinatorics. Dominique Dumont~[2], however, recently recognized 
them as the cases $k = 2,3,4$ of the following marvellous 
conjecture (which he also proved for 
$n-k = 0,8,16,24,32,40$).

\medskip

\noindent
{\sc Conjecture~(Dumont).~--~}
{\it The following relation holds for all positive integers $n$ and $k$:} 
$$
r_k^{1(2)}(n) \; = \; c_k^{1(4)}(n) - (-1)^k c_k^{3(4)}(n).
$$

\medskip

Once the right statement has been found, the proof almost takes care of itself.
We think that this lemma can be applied directly to Dumont's conjecture. 
In any case, it can be applied to our main result.

\medskip

{\sc Theorem.~--~}
{\it Let the infinite matrices $A = (a_{ij})$ and $B = (b_{ij})$, $i,j = 0,1,2,3,\dots$
be defined by $a_{ij} := q^{(4i+1)(4j+1)}$ for all $i,j$ and by $b_{ij} := -q^{(4i-1)(4j-1)}$, 
$b_{0j} = b_{i0} = 0$ for $i,j > 0$, whereas $b_{00} := \sum_{n=0}^\infty q^{(2n+1)^2}$.
Then there exists an inversible matrix $X$ such that $XB = AX$. In particular,}
$\hbox{tr}\bigl[A^k\bigr] = \hbox{tr}\bigl[B^k\bigr]$, {\it which is the generating
function formulation of the identity $c_k^{1(4)}(n) = r_k^{1(2)}(n) + (-1)^k c_k^{3(4)}(n)$.}

\bigskip


\vskip 6pt
\centerline{\hbox to 2cm{\hrulefill}}
\vskip 9pt



\bigskip

\noindent{\bf 1. Combinatoire}

\medskip

En suivant Andrews~[1], commen\c cons par rappeler les plus beaux th\'eor\`emes 
de la combinatoire des $q$-s\'eries formelles (ceux qui pr\'ef\`erent l'analyse 
complexe vont poser $|q|<1$), qui utilisent tous les notations
$$
(a;q)_n := (1-a)(1-aq)\cdots(1-aq^{n-1}), \quad (a;q)_\infty := \lim_{n\to\infty} (a;q)_n, \quad (a;q)_0 := 1.
$$ 

\eject

\noindent
{\sc Th\'eor\`eme~$q$-Binomial~(Cauchy).} \quad
$$
\sum_{n=0}^\infty {(a;q)_n \over (q;q)_n} \; t^n 
\; = \; 
{(at;q)_\infty \over (t;q)_\infty}.
$$

\medskip

{\it D\'emonstration. -- }
Posons $F(t) := (at;q)_\infty / (t;q)_\infty =: \sum_{n=0}^\infty A_n(a;q)\;t^n$. 
On v\'erifie imm\'ediate\-ment $(1-t)F(t) = (1-at)(atq;q)_\infty / (tq;q)_\infty = (1-at)F(tq)$,
d'o\`u $A_n(a;q)-A_{n-1}(a;q) = q^nA_n(a;q)-aq^{n-1}A_{n-1}(a;q)$, i.e. 
$A_n(a;q) = A_{n-1}(a;q) (1-aq^{n-1})/(1-q^n)$.\qed 

\bigskip

\noindent
{\sc Corollaire~(Euler).} \quad
$$
\sum_{n=0}^\infty {t^n \over (q;q)_n}  \; = \; {1 \over (t;q)_\infty}, \qquad
\sum_{n=0}^\infty {t^n q^{n(n-1)/2} \over (q;q)_n}  \; = \; (-t;q)_\infty. 
$$

\medskip

{\it D\'emonstration. -- }
Il suffit de poser $a=0$ o\`u bien de remplacer $a$ par $a/b$ et $t$ par $bt$ pour poser
$b=0$ et $a=-1$ ensuite.\qed

\bigskip

\noindent
{\sc Th\'eor\`eme~Triple Produit~(Jacobi).} \quad
$$
\sum_{n=-\infty}^\infty q^{n^2} z^n \; = \; (q^2;q^2)_\infty(-qz;q^2)_\infty(-q/z;q^2)_\infty
$$

\medskip

\noindent
{\it D\'emonstration. -- }
Nous utilisons trois fois le corollaire pr\'ec\'edent:
$$\openup2pt\eqalign{
(q^2;q^2)_\infty(-qz;q^2)_\infty 
\; &= \; (q^2;q^2)_\infty \sum_{n=0}^\infty {q^{n^2}z^n \over (q^2;q^2)_n}
\;  = \; \sum_{n=0}^\infty q^{n^2}z^n (q^{2n+2};q^2)_\infty \cr
\; &= \; \sum_{n=-\infty}^\infty q^{n^2}z^n (q^{2n+2};q^2)_\infty 
\;  = \; \sum_{n=-\infty}^\infty q^{n^2}z^n \sum_{m=0}^\infty {(-1)^m q^{m^2+m+2mn} \over (q^2;q^2)_m} \cr
\; &= \; \sum_{m=0}^\infty {(-1)^m q^m z^{-m} \over (q^2;q^2)_m} \sum_{n=-\infty}^\infty q^{(n+m)^2} z^{n+m}
\;  = \; \sum_{m=0}^\infty {(-q/z)^m \over (q^2;q^2)_m} \sum_{n=-\infty}^\infty q^{n^2} z^{n} \cr
\; &= \; {1 \over (-q/z;q^2)_\infty} \sum_{n=-\infty}^\infty q^{n^2} z^{n}.\qed
}$$

\bigskip

\noindent
{\sc Corollaire~(Gau\ss).} \quad
$$
\sum_{n=-\infty}^\infty (-1)^n q^{n^2} \; = \; {(q;q)_\infty \over (-q;q)_\infty}, \qquad 
\sum_{n=0}^\infty q^{n(n+1)/2} \; = \; {(q^2;q^2)_\infty \over (q;q^2)_\infty}. 
$$

\medskip

{\it D\'emonstration. -- }
L'identit\'e $1+q^n = (1-q^{2n})/(1-q^n)$ entra\^\i ne $(-q;q)_\infty = (q^2;q^2)_\infty / (q;q)_\infty = 1/(q;q^2)_\infty$.
Les cas $z = -1$ et $z = q$ du th\'eor\`eme pr\'ec\'edent impliquent donc bien les deux identit\'es
$
\sum_{n=-\infty}^\infty (-1)^n q^{n^2} = (q^2;q^2)_\infty(q;q^2)_\infty(q;q^2)_\infty = (q;q)_\infty(q;q^2)_\infty = (q;q)_\infty / (-q;q)_\infty
$
et
$
\sum_{n=0}^\infty q^{n(n+1)/2} = {1 \over 2} \sum_{n=-\infty}^\infty q^{n(n+1)/2} = {1 \over 2} (q;q)_\infty(-q;q)_\infty(-1;q)_\infty
= (q;q)_\infty(-q;q)_\infty(-q;q)_\infty = (q^2;q^2)_\infty(-q;q)_\infty = (q^2;q^2)_\infty / (q;q^2)_\infty.\qed
$

\eject


\bigskip

\noindent{\bf 2. Alg\`ebre lin\'eaire}

\medskip

Nous regardons maintenant des matrices infinies $A = (a_{ij})$, $i,j = 0,1,2,3,\dots$, dont les \'el\'ements $a_{ij}$
sont des $q$-s\'eries formelles. Appelons une telle matrice {\sl admissible} si et seulement si, pour
tout $n \in {\bb N}$, il existe un $k \in {\bb N}$ tel que $|i-j| > k$ implique deg$(a_{ij}) > n$.
L'identit\'e $I$ est admissible et le produit $AB$ de deux matrices admissibles $A$ et $B$ est bien 
d\'efini et admissible. De plus, la multiplication des matrices admissibles est associative. 
Si $A = (a_{ij})$ est une matrice admissible avec deg$(a_{ij}) > 0$ pour tout $i,j$, alors 
$(I+A)^{-1} = \sum_{k=0}^\infty (-1)^k A^k$ est \'egalement bien d\'efini et admissible. 
Appelons la matrice $A = (a_{ij})$ {\sl admise} si et seulement si, pour tout $n \in {\bb N}$, 
il existe un $k \in {\bb N}$ tel que $\max(i,j) > k$ implique deg$(a_{ij}) > n$. Pour chaque
matrice admise tr$\bigl[A\bigr]$ est bien d\'efinie. De plus, si $A$ est une matrice admise
et $X$ est une matrice admissible, alors $AX$ et $XA$ sont des matrices admises et l'on
a l'identit\'e $\hbox{tr}\bigl[AX\bigr] =  \hbox{tr}\bigl[XA\bigr]$. 

Nous nous int\'eressons ici plus sp\'ecialement aux matrices de la forme $q^{ij}$
(pour les \'etudier dans le cas $|q| < 1$ et mettre fin \`a l'injustice {\lfq fourielle\rfq}
de se borner au bord $|q| = 1$). D\'efinissons donc les deux matrices admises $A = (a_{ij})$, 
$B = (b_{ij})$ et la matrice admissible $X = (x_{ij})$ par 
$$
a_{ij} := q^{(4i+1)(4j+1)}, 
$$
$$
b_{ij} := \openup 2pt\left\{\eqalign{
          -q^{(4i-1)(4j-1)},                \;\; & \hbox{si} \;\; i,j > 0, \cr
          \sum_{n=0}^\infty q^{(2n+1)^2}, \;\;\; & \hbox{si} \;\; i=j=0, \cr
          0,             \qquad \quad         \; & \hbox{sinon,} \cr  
}\right. \qquad\qquad
x_{ij} := \openup 2pt\left\{\eqalign{
          -{q^{12(j-i-1)+4} \over 1-q^{16(j-i-1)+8}},                \;\; & \hbox{si} \;\; j > i, \cr
          {q^{4(i-j)} \over 1-q^{16(i-j)+8}},                  \quad \;\; & \hbox{si} \;\; 1 \le j \le i, \cr
          {(q^8;q^{16})_i \over (q^{16};q^{16})_i} \; q^{4i},  \quad \;\; & \hbox{si} \;\; j = 0. \cr  
}\right. 
$$

\bigskip

\noindent
{\sc Th\'eor\`eme.~--~}
{\it Nous avons $\; XB \; = \; AX$.}

\medskip

\noindent
{\it D\'emonstration. -- }
D'abord, il nous faut montrer, pour tout $i \ge 0$ et $j \ge 1$, que 
$$\openup2pt\eqalign{
     &-\sum_{m=1}^i {q^{4(i-m)} \over 1-q^{16(i-m)+8}} \cdot q^{(4m-1)(4j-1)}
      +\sum_{m=i+1}^\infty {q^{12(m-i-1)+4} \over 1-q^{16(m-i-1)+8}} \cdot q^{(4m-1)(4j-1)} \cr
= \; &-\sum_{n=0}^{j-1} {q^{12(j-n-1)+4} \over 1-q^{16(j-n-1)+8}} \cdot q^{(4n+1)(4i+1)}
      +\sum_{n=j}^\infty {q^{4(n-j)} \over 1-q^{16(n-j)+8}} \cdot q^{(4n+1)(4i+1)}. \cr
}$$
En posant $m-i-1 = n-j = k$, nous obtenons 
$$\openup2pt\eqalign{
         &\sum_{n=j}^\infty {q^{4(n-j)} \over 1-q^{16(n-j)+8}} \cdot q^{(4n+1)(4i+1)}
          -\sum_{m=i+1}^\infty {q^{12(m-i-1)+4} \over 1-q^{16(m-i-1)+8}} \cdot q^{(4m-1)(4j-1)} \cr
   = \; &q^{16ij+4j-4i+1}\sum_{k=0}^\infty q^{8k} { (q^{16k+8})^i - (q^{16k+8})^j \over 1-q^{16k+8} } 
\; = \;  q^{16ij+4j-4i+1}\sum_{k=0}^\infty q^{8k} (-1)^{[i > j]} \sum_{l=\min(i,j)}^{\max(i,j)-1} (q^{16k+8})^l \cr
   = \; &(-1)^{[i > j]} q^{16ij+4j-4i+1} \sum_{l=\min(i,j)}^{\max(i,j)-1} { q^{8l} \over 1-q^{16l+8} }, \qquad \hbox{o\`u} \;\;
[i > j] := \openup 2pt\left\{\eqalign{
           1, \;\; & \hbox{si} \;\; i > j, \cr
           0, \;\; & \hbox{sinon}. \cr
           }\right.   \cr
}$$
D'autre part, en posant $i-m = j-n-1 = k$, nous obtenons
$$\openup2pt\eqalign{
     &\sum_{m=1}^i {q^{4(i-m)} \over 1-q^{16(i-m)+8}} \cdot q^{(4m-1)(4j-1)}
      -\sum_{n=0}^{j-1} {q^{12(j-n-1)+4} \over 1-q^{16(j-n-1)+8}} \cdot q^{(4n+1)(4i+1)} \cr
= \; &q^{16ij+4j-4i+1} \Biggl[ (-1)^{[j > i]} \sum_{k=\min(i,j)}^{\max(i,j)-1} q^{8k} { (q^{16k+8})^{-\min(i,j)} \over 1-q^{16k+8} }
      +\sum_{k=0}^{\min(i,j)-1} q^{8k} { (q^{16k+8})^{-j} - (q^{16k+8})^{-i} \over 1-q^{16k+8} } \Biggr] \cr
= \; &(-1)^{[j > i]} q^{16ij+4j-4i+1} \Biggl[ \sum_{k=\min(i,j)}^{\max(i,j)-1} q^{8k} { (q^{16k+8})^{-\min(i,j)} \over 1-q^{16k+8} }
      -\sum_{k=0}^{\min(i,j)-1} q^{8k} \sum_{l=\min(i,j)}^{\max(i,j)-1} (q^{16k+8})^{-l-1} \Biggr] \cr
= \; &(-1)^{[j > i]} q^{16ij+4j-4i+1} \Biggl[ \sum_{l=\min(i,j)}^{\max(i,j)-1} q^{8l} { (q^{16l+8})^{-\min(i,j)} \over 1-q^{16l+8} }
      - \sum_{l=\min(i,j)}^{\max(i,j)-1} q^{-8(l+1)} \sum_{k=0}^{\min(i,j)-1} (q^{16l+8})^{-k} \Biggr] \cr
= \; &(-1)^{[j > i]} q^{16ij+4j-4i+1} \sum_{l=\min(i,j)}^{\max(i,j)-1} \Biggl[ q^{8l} { (q^{16l+8})^{-\min(i,j)} \over 1-q^{16l+8} }
      - q^{-8(l+1)} { (q^{16l+8})^{-\min(i,j)+1} - q^{16l+8} \over 1 - q^{16l+8} } \Biggr] \cr 
= \; &(-1)^{[j > i]} q^{16ij+4j-4i+1} \sum_{l=\min(i,j)}^{\max(i,j)-1} { q^{8l} \over 1-q^{16l+8} }, \cr
}$$
ce qui ach\`eve la d\'emonstration dans le cas $i \ge 0$ et $j \ge 1$. Dans le cas $j=0$, gr\^ace aux identit\'es de Cauchy et de Gau\ss,
nous avons pour tout~$i$
$$
        \sum_{n=0}^\infty {(q^8;q^{16})_n \over (q^{16};q^{16})_n} \; q^{4n} \cdot q^{(4n+1)(4i+1)}
\; = \; q^{4i+1} \sum_{n=0}^\infty {(q^8;q^{16})_n \over (q^{16};q^{16})_n} \; (q^{16i+8})^n
\; = \; q^{4i+1} {(q^{16i+16};q^{16})_\infty \over (q^{16i+8};q^{16})_\infty}
$$
$$
\; = \; q \; {(q^{16};q^{16})_\infty \over (q^{8};q^{16})_\infty} \cdot {(q^8;q^{16})_i \over (q^{16};q^{16})_i} \; q^{4i} 
\; = \; \Biggl[ \sum_{n=0}^\infty q^{(2n+1)^2} \Biggr] \cdot {(q^8;q^{16})_i \over (q^{16};q^{16})_i} \; q^{4i}. \qed
$$

\bigskip

\noindent
{\sc Corollaire.~--~}
{\it Nous avons $\; c_k^{1(4)}(n) \; = \; r_k^{1(2)}(n) + (-1)^k c_k^{3(4)}(n)$.}

\medskip

\noindent
{\it D\'emonstration. -- }
$\hbox{tr}\bigl[A^k\bigr] = \hbox{tr}\bigl[(XBX^{-1})^k\bigr] = \hbox{tr}\bigl[(XB^k)X^{-1}\bigr] 
= \hbox{tr}\bigl[X^{-1}(XB^k)\bigr] = \hbox{tr}\bigl[B^k\bigr]$.\qed

\bigskip


\bigskip

\noindent{\bf 3. Th\'eorie des nombres}

\medskip

Terminons cette Note avec quelques applications de notre th\'eor\`eme principal,
en commen\c cant avec le th\'eor\`eme dit {\lfq{\sl Eur\^eka}\rfq} de Gau\ss. 

\bigskip

\noindent
{\sc Corollaire~(Gau\ss).~--~}
{\it Tout nombre naturel se d\'ecompose en trois nombres triangulaires.}

\medskip

{\it D\'emonstration. -- }
Puisque $m = {y_1 \choose 2}+{y_2 \choose 2}+{y_3 \choose 2} \; \Leftrightarrow \; 8m+3 = (2y_1-1)^2+(2y_2-1)^2+(2y_3-1)^2$,
il faut montrer que $r_3^{1(2)}(8m+3) = c_3^{1(4)}(8m+3)+c_3^{3(4)}(8m+3) > 0$. En effet, en posant $x_1 = x_2 = 1$ et
$x_3 = 4m+1$, nous avons $x_1x_2+x_2x_3+x_3x_1 = 8m+3$.\qed

\bigskip

\noindent
{\sc Corollaire~(Jacobi).} \quad
$$
r_2^{1(2)}(8m+2) \; = \; \sum_{d|4m+1} (-1)^{(d-1)/2}.
$$

\medskip

\noindent
{\it D\'emonstration. -- } Nous avons bien $8m+2 = x_1x_2+x_2x_1 \; \Leftrightarrow \; 4m+1 = x_1x_2$.\qed

\bigskip

\noindent
{\sc Corollaire~(Jacobi).} \quad
$$
r_4^{1(2)}(8m+4) \; = \; \sum_{d|2m+1} d.
$$

\medskip

{\it D\'emonstration. -- } Nous avons bien 
$8m+2 = x_1x_2+x_2x_3+x_3x_4+x_4x_1 \; \Leftrightarrow \; 2m+1 = {x_1+x_3 \over 2} \cdot {x_2+x_4 \over 2}$.
Il s'ensuit que 
$r_4^{1(2)}(8m+4) = c_4^{1(4)}(8m+4)-c_4^{3(4)}(8m+4) 
= \sum_{d \cdot d' = 2m+1} {2d+2 \over 4} \cdot {2d'+2 \over 4} - {2d-2 \over 4} \cdot {2d'-2 \over 4}
= \sum_{d \cdot d' = 2m+1} {d+d' \over 2} = \sum_{d|2m+1} d$.\qed

\bigskip

{\sc Corollaire~(Jacobi).~--~}
{\it Notons $r_2(n) := |\{(x_1,x_2) \in {\bb Z}^2 \; | \; n = x_1^2+x_2^2 \}|$ pour tout $n \in {\bb N}$.
Alors}
$$
r_2(n) \; = \; 4 \sum_{d|n, \;2 \nmid d} (-1)^{(d-1)/2}.
$$

\medskip

{\it D\'emonstration. -- }
L'\'equivalence $2n = (x_1+x_2)^2 + (x_1-x_2)^2 \; \Leftrightarrow \; n = x_1^2+x_2^2$ fournit une preuve
bijective que $r_2(2n) = r_2(n)$, parce que l'on a automatiquement $y_1 \equiv y_2 \!\!\! \pmod 2$ si $2n = y_1^2+y_2^2$.
Il suffit donc de d\'emontrer le corollaire dans le cas $n \equiv 2 \!\!\! \pmod 4$. Dans ce cas, cependant, nous
avons $r_2(n) = 4r_2^{1(2)}(n)$.\qed

\bigskip

{\sc Corollaire~(Jacobi).~--~}
{\it Notons $r_4(n) := |\{(x_1,x_2,x_3,x_4) \in {\bb Z}^4 \; | \; n = x_1^2+x_2^2+x_3^2+x_4^2 \}|$
pour tout $n \in {\bb N}$. Alors}
$$
r_4(n) \; = \; 8 \sum_{d|n, \;4 \nmid d} d.
$$

\medskip

{\it D\'emonstration. -- }
L'\'equivalence $2n = (x_1+x_2)^2 + (x_1-x_2)^2 + (x_3+x_4)^2 + (x_3-x_4)^2 \; \Leftrightarrow \; n = x_1^2+x_2^2+x_3^2+x_4^2$ 
fournit une preuve bijective que $r_2(2n) = r_2(n)$ si $n$ est pair, parce que l'on a automatiquement 
$y_1 \equiv y_2 \equiv y_3 \equiv y_4 \!\!\! \pmod 2$ si $2n = y_1^2+y_2^2+y_3^2+y_4^2$.
Si $n$ est impair, alors exactement un tiers des solutions de $2n = y_1^2+y_2^2+y_3^2+y_4^2$ satisfont aux relations 
$y_1 \equiv y_2 \!\!\! \pmod 2$ et $y_3 \equiv y_4 \!\!\! \pmod 2$, c'est-\`a-dire $r_2(2n) = 3r_2(n)$.
Comme $4 \nmid d$ assure les m\^emes relations pour le membre de droite, il suffit de d\'emontrer le 
corollaire dans le cas $n \equiv 4 \!\!\! \pmod 8$. Dans ce cas, cependant, nous avons 
$r_4(n) = 16r_4^{1(2)}(n) + r_4(n/4) = 16r_4^{1(2)}(n) + r_4(n)/3$, i.e. $r_4(n) = 24r_4^{1(2)}(n)$.\qed

\bigskip

\noindent
{\sc Corollaire~(Lagrange).~--~}
{\it Tout nombre naturel se d\'ecompose en quatre carr\'es.}\qed

\bigskip

{\bf Remerciements.} 
En tout premier lieu, je voudrais remercier vivement Dominique Dumont d'avoir rendu possible la 
r\'edaction de cet article en pr\'esentant sa conjecture merveilleuse au s\'eminaire
{\lfq Th\'eorie des nombres et Combinatoire\rfq} \`a l'Institut Camille Jordan.
Je remercie aussi Victor Guo pour des remarques fort utiles.

\bigskip


\centerline{\bf R\'ef\'erences bibliographiques}

\medskip

{\baselineskip=9pt\eightrm
\item{[1]}Andrews G., The theory of partitions, Cambridge University Press, 1998;
en russe: Teoriya razbienii, Nauka, 1982. 
\smallskip 

\item{[2]}Dumont D., A Conjecture on Sums of Any Number of Odd Squares, pr\'epublication.
\smallskip 

\item{[3]}Weil A., Sur les sommes de trois et quatre carr\'es,
Enseignement Math.~II.~S\'er.~20 (1974) 215-222.
\smallskip
}
\end